\documentclass[11pt]{article}
\usepackage{amssymb,amsmath}

\newcommand{\be}{\begin{equation}}
\newcommand{\ee}{\end{equation}}
\newcommand{\bea}{\begin{eqnarray}}
\newcommand{\eea}{\end{eqnarray}}
\newtheorem{theorem}{Theorem}

\newtheorem{lemma}{Lemma}

\def\1#1{^{(#1)}}

\begin{document}
\title{On Summatory Totient Functions}
\author{Leonid G. Fel\\
Department of Civil Engineering, Technion, Haifa 32000, Israel\\
\vspace{-.2cm}
\\{\sl e-mail: lfel@tx.technion.ac.il}}
\date{\today}
\maketitle
\def\be{\begin{equation}}
\def\ee{\end{equation}}
\def\bea{\begin{eqnarray}}
\def\p{\prime}
\vspace{-1cm}
\begin{abstract}
The lower and upper bounds are found for the leading term of summatory totient 
function $\sum_{k\leq N}k^u\phi^v(k)$ in various ranges of $u\in{\mathbb R}$ 
and $v\in{\mathbb Z}$.
\\{\bf Keywords:} Summatory totient functions, Asymptotic analysis.\\
{\bf 2000 Mathematics Subject Classification:} 11N37.
\end{abstract}
{\large\bf 1.} We study the summatory totient function associated with the 
Euler function $\phi(k)$,
\begin{eqnarray}
F\left[k^u\phi^v,N\right]=\sum_{k\leq N}k^u\phi^v(k)\;,\;\;\;
u\in{\mathbb R}\;,\;\;v\in{\mathbb Z}\;.\label{er1}
\end{eqnarray}
The function $F\left[k^u\phi^v,N\right]$ has been the subject of intensive study
for the last century and is classically known \cite{apo95} for $u\leq 0,\;v=1$. 
The other results include $u=0,v=-1$ \cite{lan00}, $v=-u>0$ \cite{cho30}, 
\cite{bal96} and references therein, $v\geq 0,\;u<-v-1$ \cite{liu2}, $u=1,v=-1$ 
\cite{now89}, \cite{sit85}, $u=v=-1$ \cite{mor00}, \cite{ste76}. 
The leading and error terms for $u=0$, $v\in{\mathbb Z}_+$, were calculated in 
\cite{cho30} and \cite{bal96}, respectively. An extensive survey on the 
number-theoretical properties of $\phi(k)$ and the leading and error terms of 
some summatory functions (\ref{er1}) is presented in \cite{san06}. In this 
article we give the lower and upper bounds for the leading term of $F\left[k^u
\phi^v,N\right]$ in various ranges of $u,\;v$. 

For this purpose put the following notations,
\begin{eqnarray}
\lim_{N\to\infty}\frac{F[k^u\phi^v,N]}{N^{u+v+1}}=A(u,v)\;,\;\;\;
\lim_{N\to\infty}\frac{F[k^u\phi^v,N]}{\ln N}=B(u,v)\;,\;\;\;
\lim_{N\to\infty}F[k^u\phi^v,N]=C(u,v)\;,\label{ler1a}
\end{eqnarray}
and note that for $v=0$ these asymptotics read
\begin{eqnarray}
A(u,0)=(u+1)^{-1}\;,\;\;u>-1\;;\;\;\;\;B(u,0)=1\;,\;\;u=-1\;;\;\;\;\;   
C(u,0)=\zeta(-u)\;,\;\;u<-1\;.\nonumber
\end{eqnarray}
Here $\zeta(s)$ stands for the Riemann zeta function.

{\large\bf 2.} Start with auxiliary summatory function $F[k^uJ_v,N]=\sum_{k
\leq N}k^uJ_v(k)$ which is associated with the Jordan totient function $J_v(k)$,
\begin{eqnarray}
J_v(k)=k^v\prod_{p_j\;\mid\;k}\left(1-\frac1{p_j^v}\right)=\sum_{d\;\mid\;k}
\mu(d)\left(\frac{k}{d}\right)^v\;,\;\;\;v\in{\mathbb Z}_+\;,\;\;\;J_v\left(
p^r\right)=p^{vr}\left(1-\frac1{p^v}\right)\;,\label{jr1}
\end{eqnarray}
where $\mu(d)$ denotes the M\"obius function. The leading term of $F[k^uJ_v,N]$ 
can be calculated exactly in the different ranges $u+v>-1$, $u+v=-1$ and $u+v
<-1$. Making worth of standard analytic methods \cite{apo95} (see also 
\cite{adh90}), we get
\begin{eqnarray}
\lim_{N\to\infty}\frac{\sum_{k\leq N}k^uJ_v(k)}{N^{u+v+1}}&=&\lim_{N\to\infty}
\frac1{N^{u+v+1}}\sum_{k\leq N}k^u\sum_{d\;\mid\;k}\mu(d)\left(\frac{k}{d}
\right)^v\nonumber\\
&=&\lim_{N\to\infty}\frac1{N^{u+v+1}}\sum_{k_1\leq N}\frac{\mu(k_1)}{k_1^v}
\sum_{k_2\leq N/k_1}(k_1k_2)^{u+v}=\frac{\sum_{k_1=1}^{\infty}\mu(k_1)\cdot
k_1^{-v-1}}{u+v+1}\nonumber\\
&=&\frac1{u+v+1}\;\frac1{\zeta(v+1)}\;,\;\;\;\;\;u+v>-1\;,\label{jr2}
\end{eqnarray}
\begin{eqnarray}
\lim_{N\to\infty}\frac{\sum_{k\leq N}k^{-v-1}J_v(k)}{\ln N}&=&\lim_{N\to\infty}
\frac1{\ln N}\sum_{k\leq N}\frac1{k^{v+1}}\sum_{d\;\mid\;k}\mu(d)\left(\frac{ 
k}{d}\right)^v\nonumber\\
&=&\lim_{N\to\infty}\frac1{\ln N}\sum_{k_1\leq N}\frac{\mu(k_1)}{k_1^v}\sum_{
k_2\leq N/k_1}\frac1{k_1k_2}=\sum_{k_1=1}^{\infty}\frac{\mu(k_1)}{k_1^{v+1}}
\nonumber\\
&=&\frac1{\zeta(v+1)}\;,\;\;\;\;\;u+v=-1\;,\label{jr3}
\end{eqnarray}
\begin{eqnarray}
\lim_{N\to\infty}\sum_{k\leq N}\frac{J_v(k)}{k^{-u}}&=&\lim_{N\to\infty}\sum_{k\leq 
N}\frac1{k^{-u}}\sum_{d\;\mid\;k}\mu(d)\left(\frac{k}{d}\right)^v=\lim_{N\to
\infty}\sum_{k_1\leq N}\frac{\mu(k_1)}{k_1^v}\sum_{k_2\leq N/k_1}\frac1{
(k_1k_2)^{-(u+v)}}\nonumber\\
&=&\zeta(-u-v)\cdot \sum_{k_1=1}^{\infty}\frac{\mu(k_1)}{k_1^{-u}}=
\frac{\zeta(-u-v)}{\zeta(-u)}\;,\;\;\;\;\;u+v<-1\;.\label{jr4}
\end{eqnarray}
Hence follow the bounds for $A(u,v)$, $B(u,v)$ and $C(u,v)$ in the case 
$v\in{\mathbb Z}_+$.
\begin{lemma}\label{lem11}{\rm $\;$}\\
For $v\in{\mathbb Z}_+$ the following asymptotics hold
\begin{eqnarray}
&&\mbox{If}\;\;u+v>-1\;,\;\;\;\mbox{then}\;\;\;\;\;0<A(u,v)\leq
\frac{(u+v+1)^{-1}}{\zeta(v+1)}\;,\nonumber\\
&&\mbox{If}\;\;u+v=-1\;,\;\;\;\mbox{then}\;\;\;\;\;0<B(u,v)\leq
\frac1{\zeta(v+1)}\;,\nonumber\\
&&\mbox{If}\;\;u+v<-1\;,\;\;\;\mbox{then}\;\;\;\;\;0<C(u,v)\leq
\frac{\zeta(-u-v)}{\zeta(-u)}\;,\nonumber
\end{eqnarray}
where the upper bounds are attained iff $v=1$.
\end{lemma}
{\sf Proof} $\;\;\;$Observe that the following inequality holds
\begin{eqnarray}
\phi^v(k)=k^v\prod_{p_j\;\mid\;k}\left(1-\frac1{p_j}\right)^v\leq k^v\prod_{p_j
\;\mid\;k}\left(1-\frac1{p_j^v}\right)=J_v(k)\;,\;\;\;v\in{\mathbb Z}_+\;,
\label{ler2}
\end{eqnarray}
since $(1-x^v)/(1-x)^v=(1+x+\ldots+x^{v-1})/(1-x)^{v-1}\geq 1$ if $x<1$. The
last inequality becomes rigorous if and only if $v>1$. Combining now 
(\ref{ler2}) with (\ref{jr2}), (\ref{jr3}) and (\ref{jr4}) we arrive at the
proof of Lemma.$\;\;\;\;\;\;\Box$

Illustrate Lemma \ref{lem11} by three known examples taken from 
\cite{apo95}, p. 71,
\begin{eqnarray}
A(-\alpha,1)=\frac1{(2-\alpha)\zeta(2)},\;\;\alpha<2\;;\;\;\;\;B(-2,1)=\frac1{
\zeta(2)}\;;\;\;\;\;C(-\alpha,1)=\frac{\zeta(\alpha-1)}{\zeta(\alpha)},\;\;
\alpha>2\;.\nonumber
\end{eqnarray}  
Two other examples are taken from \cite{cho30},
\begin{eqnarray}
A(0,2)&=&\frac1{3}\prod_{p}\left(1-2\frac1{p^2}+\frac1{p^3}\right)<\frac1{3}
\prod_{p}\left(1-\frac1{p^3}\right)=\frac1{3\zeta(3)}\;,\label{ler3}\\
A(-v,v)&=&\prod_{p}\left(1-\frac1{p}\left(1-\left(1-\frac1{p}\right)^v\right)
\right)\leq\prod_{p}\left(1-\frac1{p^{v+1}}\right)=\frac1{\zeta(v+1)},\nonumber
\end{eqnarray}
where inequality becomes rigorous iff $v>1$. The last example is taken from 
\cite{liu2}, 
\begin{eqnarray}
C(-v-s,v)=\zeta(s)\prod_{p}\left(1-\frac1{p^s}\left(1-\left(1-\frac1{p}
\right)^v\right)\right)<\zeta(s)\prod_{p}\left(1-\frac1{p^{v+s}}\right)=
\frac{\zeta(s)}{\zeta(v+s)}\;,\;\;\;s>1\;.\nonumber
\end{eqnarray}
{\large\bf 3.} In the case $v\in{\mathbb Z}_-$ we represent the function 
$F[k^u\phi^v;N]$ as follows, 
\begin{eqnarray}
F[k^u\phi^v;N]=\sum_{k=1}^Nk^{u+v}\prod_{p_j\;\mid\;k}\left(1-\frac1{p_j}\right)
^{-|v|},\;\;\;\;\prod_{p_j\;\mid\;k}\left(1-\frac1{p_j}\right)^{-|v|}>1\;,\;\;
\;v\in{\mathbb Z}_-\;,\label{ler5}
\end{eqnarray}
and prove Lemma on lower bounds.
\begin{lemma}\label{lem12}{\rm $\;$}\\
For $v\in{\mathbb Z}_-$ the following asymptotics hold
\begin{eqnarray}
A(u,v)>\frac1{u+v+1},\;u+v>-1;\;\;B(u,v)>1,\;u+v=-1;\;\;C(u,v)>\zeta(-u-v),\;
u+v<-1\;.\nonumber
\end{eqnarray}
\end{lemma}
{\sf Proof} $\;\;\;$In accordance with definition (\ref{ler1a}) and inequality
({\ref{ler5}) calculate the lower bound for different signs of $u+v+1$,
\begin{eqnarray}
&&1)\;\;u+v>-1\;,\;\;\;\;\;A(u,v)>\lim_{N\to \infty}\frac1{N^{u+v+1}}\sum_{k=1}
^Nk^{u+v}=\frac1{u+v+1}\;,\nonumber\\
&&2)\;\;u+v=-1\;,\;\;\;\;\;B(u,v)>\lim_{N\to \infty}\frac1{\ln N}\sum_{k=1}^N
\frac1{k}=1\;,\nonumber\\
&&3)\;\;u+v<-1\;,\;\;\;\;\;C(u,v)>\sum_{k=1}^N\frac1{k^{-(u+v)}}=\zeta(-u-v)\;.
\nonumber
\end{eqnarray}
Lemma is proven.$\;\;\;\;\;\;\Box$

As for the upper bounds, the problem is much more difficult than in the case of
nonnegative $v$. There are different ways to find the bounds applying the 
Tauberian theorem to the corresponding Dirichlet series or making use of 
inequalities for arithmetic functions
\footnote{Based on the Tauberian theorem Z. Rudnick \cite{rud07} gave an elegant
proof of convergence of summatory function $F(u,v;N)/\ln N$, $u+v=-1$, and
calculated its leading term. As for the 2nd approach, in Section 5 we give 
another proof of convergence of summatory function $F(u,v;N)$, $u+v<-1$, based 
on two inequalities for the Euler totient $\phi(k)$ and divisor $\sigma(k)$ 
functions.}.
In this article we follow the refined proof of the Landau theorem \cite{lan00} 
given in \cite{kon80}.
\begin{lemma}\label{lem13}{\rm $\;$}\\
For $v\in{\mathbb Z}_-$ the following asymptotics hold
\begin{eqnarray}
&&\mbox{If}\;\;u+v>-1\;,\;\;\;\mbox{then}\;\;\;\;\;\;\;A(u,v)<2^{\frac{|v|}{2}}
\cdot {\cal D}_{\infty}(v,1)\cdot (u+v+1)^{-1}\;,\nonumber\\   
&&\mbox{If}\;\;u+v=-1\;,\;\;\;\mbox{then}\;\;\;\;\;\;\;\;B(u,v)<2^{\frac{|v|}
{2}}\cdot {\cal D}_{\infty}(v,1)\;,\nonumber\\
&&\mbox{If}\;\;u+v<-1\;,\;\;\;\mbox{then}\;\;\;\;\;\;\;\;C(u,v)<2^{\frac{|v|}
{2}}\cdot {\cal D}_{\infty}(v,-u-v)\cdot\zeta(-u-v)\;,\nonumber   
\end{eqnarray}
where ${\cal D}_{\infty}(v,s)=\prod_{r=1}^{|v|}\zeta\left(s+r/2\right)$.
\end{lemma}
{\sf Proof} $\;\;\;$Consider a summatory function
\begin{eqnarray}
F\left[\frac{f(k)}{\phi^m(k)};N\right]=\sum_{k\leq N}\frac{f(k)}{\phi^m(k)}\;,
\label{l6a}
\end{eqnarray}
where $f(k)$ is completely multiplicative function. Notice \cite{apo95} that
\begin{eqnarray}
\frac{k}{\phi(k)}=\sum_{d\;\mid\;k}\frac{\mu^2(d)}{\phi(d)}\;,\label{ko22}
\end{eqnarray}
where a sum is taken over all divisors $d$ of $k$. Make use of (\ref{ko22}) in 
summation identity \cite{kon80}
\begin{eqnarray}
F\left[\frac{f(k)}{\phi^m(k)};N\right]&=&\sum_{k_1\leq N}\frac{f(k_1)}
{k_1\phi^{m-1}(k_1)}\sum_{d\;\mid\;k_1}\frac{\mu^2(d)}{\phi(d)}=\sum_{k_1k_2
\leq 
N}\frac{\mu^2(k_1)}{\phi(k_1)}\cdot\frac{f(k_1k_2)}{k_1k_2\phi^{m-1}(k_1k_2)}
\nonumber\\
&=&\sum_{k_1\leq N}\frac{\mu^2(k_1)}{\phi(k_1)}\sum_{k_2\leq N/k_1}\frac{
f(k_1k_2)}{k_1k_2\phi^{m-1}(k_1k_2)}\;,\nonumber
\end{eqnarray}
and perform a multiple summation in the last equality $m$ times
\begin{eqnarray}
F\left[\frac{f(k)}{\phi^m(k)};N\right]&=&\sum_{k_1\leq N}\frac{\mu^{2}(k_1)}
{\phi(k_1)}\left\{\sum_{k_1k_2\leq N}\frac{\mu^{2}(k_1k_2)}{\phi(k_1k_2)}
\left\{\sum_{k_1k_2k_3\leq N}\frac{\mu^{2}(k_1k_2k_3)}{\phi(k_1k_2k_3)}
\;\ldots\;\right.\right.\label{l8}\\
&&\left.\left.\left\{\sum_{\prod_{i=1}^{m}k_i\leq N}\frac{\mu^{2}\left\{\prod_{
i=1}^{m}k_i\right)}{\phi\left(\prod_{i=1}^{m}k_i\right)}\left\{\sum_{k_{m+1}
\leq N/\Pi_m}\frac{f\left(\prod_{i=1}^{m+1}k_i\right)}{\prod_{i=1}^{m+1}k_i^{
m}}\right\}\right\}\;\ldots\;\right\}\right\}\;\;\;\;\nonumber 
\end{eqnarray}  
where $\Pi_m=\prod_{i=1}^mk_i$. Denote by $\Theta\left[k^u\phi^v(k),N\right]$ 
the last sum in (\ref{l8}) and consider it for $f(k)=k^u$ and $m=-v$, 
$v\in{\mathbb Z}_-$,
\begin{eqnarray}
\Theta\left[k^u\phi^v(k),N\right]=\sum_{k_{|v|+1}\leq N/\Pi_{|v|}}\frac{
\prod_{i=1}^{|v|+1}k_i^u}{\prod_{i=1}^{|v|+1}k_i^{-v}}=\prod_{i=1}^{|v|}
k_i^{u+v}\sum_{k_{|v|+1}\leq N/\Pi_{|v|}}k_{|v|+1}^{u+v}\;.\nonumber
\end{eqnarray}
Thus, for different signs of $u+v+1$ we have
\begin{eqnarray}
\Theta\left[k^u\phi^v(k),N\right]&=&\left(\Pi_{|v|}\right)^{u+v}\sum_{k_{|v|  
+1}\leq N/\Pi_{|v|}}k_{|v|+1}^{u+v}<\frac{N^{u+v+1}}{u+v+1}\;\frac1{\Pi_{|v|}}
\;\;\;,\;\;\mbox{if}\;\;\;\;\;u+v>-1\;,\;\;\label{l22}\\
\Theta\left[k^u\phi^v(k),N\right]&=&\frac1{\Pi_{|v|}}\sum_{k_{m+1}\leq N/
\Pi_{|v|}}\frac1{k_{|v|+1}}<\left(\ln \frac{N}{\Pi_{|v|}}+\gamma\right)\;
\frac1{\Pi_{|v|}}\;,\;\;\mbox{if}\;\;\;\;\;u+v=-1\;,\;\;\label{l23}\\
\Theta\left[k^u\phi^v(k),N\right]&=&\frac1{\left(\Pi_{|v|}\right)^{|u+v|}}
\sum_{k_{|v|+1}\leq N/\Pi_{|v|}}\frac1{k_{|v|+1}^{-u-v}}<\frac{\zeta(-u-v)}
{\left(\Pi_{|v|}\right)^{|u+v|}}\;\;\;\;\;,\;\;\mbox{if}\;\;\;\;\;u+v<-1\;.\;\;
\label{l24}
\end{eqnarray}
where $\gamma$ is the Euler-Mascheroni constant. Denote by $D_N(v,s)$ the
multiple sum
\begin{eqnarray}
D_N(v,s)=\sum_{k_1\leq N}\frac1{k_1^s}\frac{\mu^{2}(k_1)}{\phi(k_1)}\left\{
\sum_{k_1k_2\leq N}\frac1{k_2^s}\frac{\mu^{2}(k_1k_2)}{\phi(k_1k_2)}\;\ldots\;
\left\{\sum_{\prod_{i=1}^{m}k_i\leq N}\frac1{k_{|v|}^s}\frac{\mu^{2}\left\{
\prod_{i=1}^{|v|}k_i\right)}{\phi\left(\prod_{i=1}^{|v|}k_i\right)}\right\}\;
\ldots\;\right\}.\;\;\label{l25}
\end{eqnarray}
Substitute (\ref{l22}), (\ref{l23}) and (\ref{l24}) into (\ref{l8}) and take in 
mind (\ref{l25}). Thus, we get
\begin{eqnarray}
F\left[k^u\phi^v(k);N\right]<\left\{\begin{array}{lll}
 D_N(v,1)\cdot (u+v+1)^{-1}\cdot N^{u+v+1}&,&u>-v-1\;,\\
 D_N(v,1)\cdot \ln N&,& u=-v-1\;,\\
 D_N(v,-u-v)\cdot \zeta(-u-v)&,&u<-v-1\;.\end{array}\right.\label{l26}
\end{eqnarray}
Consider the function $D_N(v,s)$ and make worth of elementary inequalities for 
the M\"obius function $\mu^{2}(k)\leq 1$ and for the Euler function \cite{osb65}
\begin{eqnarray}
\phi(k)\geq \sqrt{k}\;,\;\;\;\;\mbox{if}\;\;\;\;k\neq 2,\;6\;.\label{os1}
\end{eqnarray}  
There are two ways how to exploit (\ref{os1}) in order to get the upper bound
for $D_N(v,s)$. One of them is to calculate two separate terms for $k=2$ and
$k=6$ in every sum of (\ref{l25}) and to apply $\phi(k)\geq\sqrt{k}$ to the
rest of the terms. This way can provide with very tight bounds, however it needs
a lot of arithmetics and gives cumbersome formulas (see Section 4). More 
sympathetic is a way to make (\ref{os1}) less strong but more universal
\begin{eqnarray}
\sqrt{2}\cdot \phi(k)\geq \sqrt{k}\;,\;\;\;\;k\geq 1\;.\label{os2}
\end{eqnarray}
This leads to the simple expression of the bounds and is sufficient to prove a
convergence of the multiple sum in (\ref{l25}). Indeed, we have
\begin{eqnarray}
D_N(v,s)&<&\sum_{k_1\leq N}\frac{\sqrt{2}}{k_1^s\sqrt{k_1}}\left\{\sum_{k_1k_2
\leq N}\frac{\sqrt{2}}{k_2^s\sqrt{k_1k_2}}\left\{\;\ldots\;\left\{\sum_{\prod_{
i=1}^{|v|}k_i\leq N}\frac{\sqrt{2}}{k_{|v|}^s\prod_{i=1}^{|v|}\sqrt{k_i}}
\right\}\;\ldots\;\right\}\right\}\nonumber\\
&<&2^{\frac{|v|}{2}}\sum_{k_1=1}^Nk_1^{-(s+\frac{|v|}{2})}\cdot \sum_{k_2=1}^N
k_2^{-(s+\frac{|v|-1}{2})}\cdot\;\ldots\;\cdot \sum_{k_{|v|}=1}^Nk_{|v|}^{-(s+
\frac{1}{2})}<2^{\frac{|v|}{2}}\cdot {\cal D}_{\infty}(v,s)\;,\;\;\;\;\;\;\;\;
\label{l27}
\end{eqnarray}
where
\begin{eqnarray}
{\cal D}_{\infty}(v,s)=\sum_{k_1=1}^{\infty}k_1^{-(s+\frac{|v|}{2})}\cdot
\sum_{k_2=1}^{\infty}k_2^{-(s+\frac{|v|-1}{2})}\cdot\;\ldots\;\cdot
\sum_{k_{|v|}=1}^{\infty}k_{|v|}^{-(s+\frac{1}{2})}=\prod_{r=1}^{|v|}\zeta
\left(s+\frac{r}{2}\right)\;.\label{l28}
\end{eqnarray}
Combining now (\ref{l27}) and (\ref{l26}) and taking the limit $N\to\infty$ in
the latter we arrive at the upper bounds for any value of $u+v+1$.
$\;\;\;\;\;\;\Box$

We illustrate Lemma \ref{lem13} by three known examples taken from \cite{now89},
\cite{lan00} and \cite{slo}, seq. A065483, respectively,
\begin{eqnarray}
A(1,-1)=\frac{\zeta(2)\zeta(3)}{\zeta(6)}\;,\;\;\;\;B(0,-1)=\frac{\zeta(2)
\zeta(3)}{\zeta(6)}\;,\;\;\;\;C(-1,-1)=g\;\zeta(2)\;,\label{l29}
\end{eqnarray}
where $g=\prod_{p}\left[1+p^{-2}(p-1)^{-1}\right]\simeq 1.3398$ and $\zeta(2) 
\zeta(3)/\zeta(6)\simeq 1.9436$. All three constants satisfy quite well Lemma
\ref{lem13},
\begin{eqnarray}
1.9436<\sqrt{2}{\cal D}_{\infty}(-1,1)=\sqrt{2}\zeta\left(\frac{3}{2}\right)
\simeq 3.694,\;\;\;\;1.3398<\sqrt{2}{\cal D}_{\infty}(-1,2)=\sqrt{2}\zeta
\left(\frac{5}{2}\right)\simeq 1.897.\;\;\label{l40}
\end{eqnarray}
{\large\bf 4.} In this Section we derive the upper bound for $D_N(v,s)$ defined 
in (\ref{l25}) in the case $v=-1$ and show that one can improve (\ref{l27}) 
significantly. Indeed, we have
\begin{eqnarray}
D_N(-1,s)=\sum_{k\leq N}\frac1{k^s}\frac{\mu^{2}(k)}{\phi(k)}<
\sum_{k\leq N}\frac1{k^s\phi(k)}=\frac1{2^s}+\frac1{2\cdot 6^s}+
\sum_{k\leq N\atop k\neq 2,6}\frac1{k^s\phi(k)}\;.\label{ap1}
\end{eqnarray}  
Applying inequality (\ref{os1}) to the last sum in (\ref{ap1}) we get
\begin{eqnarray}
D_N(-1,s)<\frac1{2^s}+\frac1{2\cdot 6^s}+\sum_{k\leq N\atop k\neq 2,6}k^{-(s+   
\frac1{2})}=\frac1{2^s}\left(1-\frac1{\sqrt{2}}\right)+\frac1{6^s}\left(\frac1{
2}-\frac1{\sqrt{6}}\right)+\zeta\left(s+\frac1{2}\right)\;.\label{ap2}
\end{eqnarray}  
One can verify that the upper bound (\ref{ap2}) is stronger than $\sqrt{2}\;
\zeta\left(s+\frac1{2}\right)$ which follows by (\ref{l27}). Indeed, return to
(\ref{l29}) and write new upper bounds in accordance with (\ref{ap2}),
\begin{eqnarray}
&&1.9436<\frac1{2}\left(1-\frac1{\sqrt{2}}\right)+\frac1{6}\left(\frac1{2}-
\frac1{\sqrt{6}}\right)+\zeta\left(1+\frac1{2}\right)=2.774\;,\nonumber\\
&&1.3398<\frac1{2^2}\left(1-\frac1{\sqrt{2}}\right)+\frac1{6^2}\left(\frac1{2}- 
\frac1{\sqrt{6}}\right)+\zeta\left(2+\frac1{2}\right)=1.417\;,\label{ap3}
\end{eqnarray}
that is much better then $3.694$ and $1.897$ found in (\ref{l40}).

However, further evaluation of the upper bounds in the case $v<-1$ leads to   
extremely long and sophisticated formulas which always can be calculated for any
given negative integer $v$.

{\large\bf 5.} In this Section we give the upper bound for the summatory 
function $\sum_{k\leq N}k^u\phi^v(k)$, $v<0$, $u+v<-1$, making worth of the 
Robin's theorem \cite{rob84} for the divisor function $\sigma(k)$. 
\begin{theorem}\label{the19}{\rm $\;$}\\
If $v< 0$ and $u+v<-1$ then
\begin{eqnarray}
C(u,v)<E_m(u,v,\eta)+e^{\gamma |v|}\cdot\zeta^{|v|}(2)\cdot\sum_{r=0}^{|v|}
(-1)^r{|v|\choose r}\eta^{|v|-r}\frac{d^r\zeta(s)}{ds^r}_{s=-u-v}\;,\label{s1f9}
\end{eqnarray}
where $\eta=2.8651$ and
\begin{eqnarray} 
E_m(u,v,\eta)=\sum_{k=1}^{m-1}\left(k^u\phi^v(k)-e^{\gamma |v|}\cdot\zeta^{|v|}
(2)\cdot\frac{(\eta+\ln k)^{|v|}}{k^{-u-v}}\right)\;,\;\;\;m\geq 3\;.\nonumber
\end{eqnarray}
\end{theorem}
{\sf Proof} $\;\;\;$Start with known inequality \cite{apo95}
\begin{eqnarray}
\frac{k^2}{\zeta(2)}<\phi(k)\sigma(k)<k^2\;,\label{ro1}
\end{eqnarray}  
where $\sigma(k)$ denotes the divisor function and satisfies the Robin's
theorem \cite{rob84}
\begin{eqnarray}
\frac{\sigma(k)}{k}<e^{\gamma}\ln\ln k+\frac{D}{\ln\ln k}\;,\;\;\;
k\geq 3\;,\;\;\;D=0.6482\ldots\;\;.\label{ro2}
\end{eqnarray}
Making use of elementary inequalities
\begin{eqnarray}
0<\ln\ln k-\ln\ln 3<\ln k-\ln 3\;,\;\;\;\;\;k> 3\;,\nonumber
\end{eqnarray}
we combine both inequalities (\ref{ro1}) and (\ref{ro2}) which give together
\footnote{There is another similar inequality \cite{ros62}, $k/\phi(k)<e^{
\gamma}\ln\ln k+2.50637/\ln\ln k$, $k\geq 3$, which can be used for estimation 
of $C(u,v)$ by the same procedure with a similar precision.}
\begin{eqnarray}
\frac{e^{-\gamma}}{\zeta(2)}\frac1{\phi(k)}<\frac{\ln\ln k}{k}+\frac{De^{-
\gamma}}{k\ln\ln k}<\frac{\ln k-\beta}{k}+\frac{De^{-\gamma}}{k\ln\ln 3}=
\frac{\ln k+\eta}{k}\;,\label{ro3}
\end{eqnarray}
where $\beta=\ln 3-\ln\ln 3=1.00456$ and $\eta=De^{-\gamma}/\ln\ln 3-\beta=
2.8651$. Then we have
\begin{eqnarray}
\lim_{N\to\infty}F(k^u\phi^v;N)&<&\sum_{k=1}^{m-1}k^u\phi^v(k)+\lim_{N\to
\infty}e^{\gamma |v|}\;\zeta^{|v|}(2)\sum_{k=m}^N\frac{(\eta+\ln k)^{|v|}}
{k^{-u-v}}\nonumber\\
&=&E_m(u,v,\eta)+e^{\gamma |v|}\;\zeta^{|v|}(2)\lim_{N\to\infty}
\sum_{k=1}^N\frac{(\eta+\ln k)^{|v|}}{k^{-u-v}}\;,\label{ro4}
\end{eqnarray}
where $m\geq 3$ and $E_m(u,v,\eta)$ is given by
\begin{eqnarray}
E_m(u,v,\eta)=\sum_{k=1}^{m-1}k^u\phi^v(k)-e^{\gamma |v|}\;\zeta^{|v|}(2)
\sum_{k=1}^{m-1}\frac{(\eta+\ln k)^{|v|}}{k^{-u-v}}\;.\label{ro5}
\end{eqnarray}
Consider the sum in (\ref{ro4}),
\begin{eqnarray}
\lim_{N\to\infty}\sum_{k=1}^N\sum_{r=0}^{|v|}{|v|\choose r}\frac{(\ln k)^{r}
\eta^{|v|-r}}{k^{-u-v}}=\sum_{r=0}^{|v|}{|v|\choose r}\eta^{|v|-r}\lim_{N\to
\infty}\sum_{k=1}^N\frac{(\ln k)^{r}}{k^{-u-v}}\;,\nonumber
\end{eqnarray}
and make use of the $r$-th derivative of the Riemann zeta function for $Re[s]>1$
given by
\begin{eqnarray}
\frac{d^r\zeta(s)}{ds^r}_{s=w}=(-1)^r\sum_{k=1}^{\infty}\frac{(\ln k)^{r}}
{k^w}\;,\;\;\;\;\;\;\frac{d^0\zeta(s)}{ds^0}_{s=w}=\zeta(w)\;.\nonumber
\end{eqnarray}
Thus, we get
\begin{eqnarray}
C(u,v)<E_m(u,v,\eta)+e^{\gamma |v|}\;\zeta^{|v|}(2)\sum_{r=0}^{|v|}(-1)^r{
|v|\choose r}\eta^{|v|-r}\;\frac{d^r \zeta(s)}{d s^r}_{s=-u-v}\;,\label{ro6}
\end{eqnarray}
that proves Theorem.$\;\;\;\;\;\;\Box$

In the case $u=v=-1$ we have by Theorem \ref{the19}
\begin{eqnarray}
C(-1,-1)<E_m(-1,-1,\eta)+e^{\gamma}\;\zeta(2)\left(\eta\;\zeta(2)-
\zeta^{\p}(2)\right)\;,\label{ro7}
\end{eqnarray}
where according to \cite{slo}, seq. A073002, the derivative $\zeta^{\p}(2)$ is 
given by
\begin{eqnarray}
\zeta^{\p}(2)=\zeta(2)\cdot\left(\gamma+\ln (2\pi)-12\ln A_{GK}\right)=
-0.937548\;,\label{ro8}
\end{eqnarray}
and $A_{GK}=1.282427$ stands for the Glaisher-Kinkelin constant 
\cite{slo}, seq. A074962.
 
Keeping in mind (\ref{l29}) and (\ref{ro8}) rewrite (\ref{ro7}) in the form
\begin{eqnarray}
g<\frac{E_m(-1,-1,\eta)}{\zeta(2)}+10.064\;,\label{ro9}
\end{eqnarray}
and compare this upper bound with (\ref{l40}) and (\ref{ap3}). The numerical 
calculations show that (\ref{ro9}) is stronger than (\ref{l40}) and (\ref{ap3}) 
for $m\geq 20$ and $m\geq 195$, respectively.
\section*{Acknowledgement}
The useful discussion with Z. Rudnick is highly appreciated.

\end{document}